\documentclass{article}

\def\eps{\varepsilon}

\def \N {{\mathbf {N}}}

\title{ \bf Медленное убывание корреляций \\ для типичных  перемешивающих   автоморфизмов}
\author{ В.В. Рыжиков}
\date{}

\usepackage[T1]{fontenc} 
\usepackage[cp1251]{inputenc}  
\usepackage[russian]{babel}
\begin{document}

	\maketitle
	
	\begin{abstract}   
		\bf Slow decay of correlations  for generic mixing automorphisms. \rm Given $\psi(n)\to +0$ and nonzero $f\in L_2(\mu)$,  we prove 
		that for the  Alpern-Tikhonov  generic mixing automorphism $T$ of $(X,\mu)$ the set $\{n\,:\, |(T^nf,f)|>\psi(n)\}$ is infinite.
		
		Для произвольной   последовательности $\psi(n)\to +0$  и интегрируемой в квадрате  ненулевой функци $f$ для типичного относительно метрики  Альперна-Тихонова перемешивающего автоморфизма $T$  установлено медленное убывание корреляций относительно $\psi(n)$: неравенство  $|(T^nf,f)|>\psi(n)$ выполняется для бесконечного множества значений $n$.  
		
		\end{abstract} 

\large
\section{Введение}
Пусть $T$ -- автоморфизм вероятностного пространства $(X,\mu)$,  мы также обозначаем  через $T$ отвечающий ему унитарный  оператор в $L_2(\mu)$.
В работах \cite{A},\cite{T} С. Альперн и С.В. Тихонов ввели полную метрику на множестве $Mix$ всех перемешивающих автоморфизмов  пространства $(X,\mu)$. В связи с этим   возникают задачи о типичности в смысле  Бэра тех или иных свойств перемешивающих преобразований. 
Традиционно  говорят "типичный автоморфизм обладает сингулярным спектром" (это доказано в \cite{T}) вместо  "все элементы  некоторого плотного $G_\delta$-множества обладают сингулярным спектром".

В заметке дан ответ на вопрос   И.В. Подвигина  о медленном убывании корреляций  для типичного перемешивающего автоморфизма.  Решение использует конструкции ранга один и метод аппроксимации
перемешивающего преобразования неперемешивающими, изложенный в  \cite{07}. Свойство неперемешивающих преобразований (в нашем случае медленное убывание  корреляций) наследуются предельным перемешивающим автоморфизмом.  Теорема Баштанова \cite{B} о плотности класса сопряженности в $Mix$  позволяет из одного такого специального примера стандартным способом  получить плотное $G_\delta$-множество перемешивающих автоморфизмов с медленным убыванием  корреляций.

В заключичительной части заметки мы покажем, что гипотеза о перемешивающих преобразованиях из работы \cite{P} неверна.
Не только перемешивающие автоморфизмы, как это предполагала гипотеза, но все   мягко перемешивающие автоморфизмы (и только  они) 
 не допускают сходимость  ненулевых  средних $\frac 1 {{k_n}} \sum_{i=1}^{{k_n}}T^if (x)$  со скоростью  $o(\frac 1 {k_n})$.

\bf Замечание.  \rm Отметим, что типичный  перемешивающий 
автоморфизм обладает  рангом один \cite{B}.  Автор предложил в \cite{20}  другой метод доказательства  этого факта и получил   аналогичный  результат  для  перемешивающих потоков (\cite{20}, теорема 7.1). 
Из этого  вытекает, что такие потоки обладают тривиальным централизатором. 
  С учетом результатов работы  \cite{07} можно доказать, что  для типичных перемешивающих автоморфизмов  и потоков  их симметрические тензорные степени имеют однократный спектр. 

\section{Убывание корреляций для типичных перемешивающих автоморфизмов}

В силу эргодической теоремы следующее утверждение  содержательно для функций $f$ с нулевым средним, когда корреляции стремятся к 0.

\vspace{2mm} 
\bf Теорема 1. \it Пусть  $0\neq f\in L_2(\mu)$ и $\psi(n)\to +0$ при 
$n\to\infty$. Для  типичного  перемешивающего автоморфизма $T$   множество $\{n\,:\, |(T^nf,f)|> \psi(n)\}$  бесконечно.
 \rm

\vspace{2mm}
Отметим, что аналогичный  результат справедлив для  перемешивающих потоков.  Напомним, как определяется метрика $r$  на Mix, пространстве перемешивающих автоморфизмов.  
На группе  $Aut(\mu)$ всех атоморфизмов задана метрика Халмоша $\rho$:
$$ \rho(S,T)=\sum_i 2^{-i}\left(\mu(SA_i\Delta TA_i)+\mu(S^{-1}A_i\Delta T^{-1}A_i)\right),$$
где $\{A_i\}$  -- некоторое  фиксированное семейство,  плотное в алгебре всех $\mu$-измеримых множеств.
Введем  на $Aut(\mu)$ метрику  $d_w$:
$$ d_w(S,T)=\sum_{i,j} 2^{-i-j}\left|\mu(SA_i\cap A_j)-\mu(TA_i\cap A_j)\right|.$$
На Mix  метрика $r$ задается следующим образом:
$$  r(S,T)= \rho(S,T) + \sup_{n>0}d_w(S^n,T^n).$$
В \cite{T} доказано, что пространство $(Mix,r)$ является полным и сепарабельным.

Чтобы найти типичное множество автоморфизмов с медленным убыванием корреляций, нам нужно  предъявить пример такого автоморфизма. Говорим, что $T$ обладает \it  $(1-\eps)$-перемешиванием, \rm 
если всякий слабый предел  его степеней $T^{n_i}$, $n_i\to\infty$, ограниченный на пространство функций с нулевым средним, имеет норму, не превосходящую $\eps$.
Отметим, что 1-перемешивание совпадает с  обычным   перемешиванием.
 Автоморфизм  $T$ обладает \it$\eps$-жесткостью, \rm если для некоторой последовательности $m_k\to\infty$ выпоняется  $T^{m_k}\to \eps I$. \rm 
Для автоморфизма $T$ при одновременном выполнении  свойств $\eps$-жесткости
и $(1-\eps)$-перемешивания и функции 
$f$ с нулевым средним имеем  $$(T^{m_k}f,f)\to \eps(f,f).$$

Воспользовавшись методом работы  \cite{07}, мы предъявим последовательность автоморфизмов $T_k$ таких, что $T_k$ обладает $(1-\eps_k)$-перемешиванием и $\eps_k$-жесткостью, $\eps_k\to +0$. Для заранее фиксированного набора функций $f_i\neq 0$, $i=1,2,\dots,k$ будет выполнено
$(T_k^{m_{k}}f_i,f_i)\ > \ \psi({m_{k}}).$
Автоморфизмы   $T_k$ в смысле работы \cite{07}  аппроксимируют перемешивающий  автоморфизм $T$, который наследует  корреляции:
$$(T^{m_{k}}f_i,f_i)\ > \ \psi({m_{k}}), \ i=1,2,\dots,k. $$

Сформулируем вспомогательные утверждения и выведем из них теорему.

\vspace{2mm}
\bf Лемма 2. \it Пусть  $\psi(n)\to +0$ и задана последовательность функций   $f_k\in L_2(\mu)$, $\|f_k\|=1$.  Для  некоторого перемешивающего автоморфизма 
$T$ верно, что  для   всяких $k$, $m$  найдется $n>m$, для которого  
$$|(T^nf_i,f_i)|> \psi(n), \ i=1,2,\dots,k.$$ \rm

\vspace{2mm}
Лемму 2 докажем позже. Для доказательства теоремы нам  понадобится ее частный случай, когда $f_k=J_kf$ для некоторых автоморфизмов $J_k$. 

\vspace{2mm}
\bf Лемма 3. \it Пусть  $\psi(n)\to +0$,  $\|f\|_2=1$ и задан счетный  набор  автоморфизмов $J_i$.  Найдется перемешивающий автоморфизм $T$ такой, что   для любых
 $ k,m$   для некоторого  $n>m$   выполнено 
$$|(J_i^{-1}T^nJ_if,f)|> \psi(n),\ i=1,2,\dots,k.$$
\rm

\vspace{2mm} 
\bf Доказательство теоремы 1. \rm 
Определим открытые в Mix множества 
$$Y_{m}=\{T\in Mix:  |(T^mf,f)|> \psi(m)\}$$ 
и $G_\delta$-множество 
$$ G=\bigcap_M  \bigcup_{m>M} Y_{m}.$$
Нам нужно доказать, что $G$ плотно относительно метрики $r$.
Пусть семейство $\{J_k\}$ плотно в $Aut(\mu)$ относительно метрики 
$\rho$.  Рассмотрим  перемешивающий автоморфизм $T$ из леммы 3,
тогда  $\{J_k^{-1}TJ_k\}\subset G$.   Теорема Баштанова \cite{B} 
утверждает, что класс сопряженности $\{J^{-1}TJ\,:\, J\in Aut\}$ плотен в Mix. Но в нем плотно относительно метрики $r$ семейство $\{J_k^{-1}TJ_k: k\in\N\}$, что  вытекает, например,  из леммы 10 \cite{T}.  Поясним этот факт для полноты изложения.   
Из сходимости $J_{k'}\to_\rho J$ (здесь $k'$ -- некоторая последовательность) в силу неравенства 
$$|\mu(T^nJ_{k'}A_i\cap J_{k'}A_j)- \mu(T^nJA_i\cap JA_j)|\leq \mu (J_{k'}A_i\Delta J A_i) + \mu (J_{k'}A_j\Delta J A_j)$$
вытекает равномерная  по $n$ сходимость
$$\mu(T^nJ_{k'}A_i\cap J_{k'}A_j)\to \mu(T^nJA_i\cap JA_j).$$
 Получаем, что  $J^{-1}TJ$ в метрике $d_w$, следовательно, и в метрике $r$   приближается последовательностью $J_{k'}^{-1}TJ_{k'}$.
Значит,   семейство $\{J_k^{-1}TJ_k\}$ плотно в Mix.  Множество $G$ типично в Mix.  Теорема доказана.

\section  {Конструкции  ранга один}
Для удобства читателя мы напомним определение преобразований ранга один.  Из них мы позже выберем нужные примеры перемешивающих автоморфизмов $T$. 

Фиксируем натуральное  число $h_1\geq 1$ (высота башни на этапе $j=1$), последовательность  $r_j\to\infty$ (число колонн, на которые виртуально разрезается башня этапа $j$)   и последовательность целочисленных векторов (параметров надстроек)   
$$ \bar s_j=(s_j(1), s_j(2),\dots, s_j(r_j-1),s_j(r_j)).$$  
На шаге $j=1$ задан  полуинтервал $C_1$. Пусть на  шаге $j$  определена  
 система   непересекающихся полуинтервалов 
$C_j, TC_j,\dots, T^{h_j-1}C_j,$
причем на $C_j,  \dots, T^{h_j-2}C_j$
пребразование $T$ является параллельным переносом. Такой набор   полуинтервалов  называется башней этапа $j$, их объединение обозначается через $X_j$ и тоже называется башней.

Представим   $C_j$ в виде   дизъюнктного объединения  полуинтервалов 
$C_j^i$, $i=1,2,\dots, {r_j},$ одинаковой длины.  
Для  каждого $i=1,2,\dots, r_j$ рассмотрим так называемую колонну  
$$C_j^i, TC_j^i ,T^2 C_j^i,\dots, T^{h_j-1}C_j^i.$$
К каждой  колонне с номером $i$  добавим  $s_j(i)$  непересекающихся полуинтервалов (этажей)  длины, равной длине интервала $C_j^i$.
Полученные  наборы интервалов  при  фиксированных $i$,$j$ называем надстроенными колоннами  $X_{i,j}$. Отметим, что при фиксированном $j$ по построению колонны $X_{i,j}$   не пересекаются. Используя параллельный перенос интервалов, преобразование $T$ теперь  доопределим так, чтобы колонны $X_{i,j}$ имели вид 
$$C_j^i, TC_j^i, \dots,  T^{h_j}C_j^i, T^{h_j+1}C_j^i, \dots, T^{h_j+s_j(i)-1}C_j^i,$$
  а   верхние этажи   колонн  $X_{i,j}$ ($i<r_j$) преобразование $T$ параллельным переносом  отображало в нижние
этажи колонн $X_{i+1,j}$: 
$$T^{h_j+s_j(i)}C_j^i = C_j^{i+1}, \ 0<i<r_j.$$ 
Положив $C_{j+1}= C^1_j$, замечаем, что все указанные этажи надстроенных колонн в новых обозначениях имеют вид 
$$C_{j+1}, TC_{j+1}, T^2 C_{j+1},\dots, T^{h_{j+1}-1}C_{j+1},$$
 образуя башню  этапа $j+1$ высоты  
 $$ h_{j+1} =h_jr_j +\sum_{i=1}^{r_j}s_j(i).$$

Частичное определение преобразования $T$ 
на этапе $j$ сохраняется на всех последующих этапах. 
В результате   получаем пространство  $X=\cup_j X_j$ 
и обратимое преобразование $T:X\to X$, сохраняющее  
стандартную меру Лебега $\mu$  на $X$.

\section{Доказательство леммы 2. Вынуждение перемешивания} 
Известно, что в случае $r_j\to\infty$, $s_j(i)=i$,
соответствующая лестничная конструкция $S$  является перемешивающей
 \cite{Ad}, \cite{20}.
Пусть конструкция $T$   похожа на  $S$:
$$s_j(i)=i, \ \  
1\leq i\leq [(1-\eps_j)r_j], 0< \eps_j<1. $$  
Для  $i> [(1-\eps_j)r_j]$ параметры  $s_j(i)$ выбираем произвольно, но будем считать, что высоты $h_j$ у конструкций $S$ и  $T$ одинаковы,
т.е. суммы $\sum_{i=1}^{r_j}s_j(i)$ одинаковы для $S$ и  $T$.

 Для нас важно, что при $\eps_j\to +0$ конструкция $T$ перемешивает,
это доказывается  так же, как и для лестничной конструкции $S$.  

\vspace{3mm}
\bf  Частичная жесткость, частичное перемешивание. \rm

\vspace{3mm}
Пусть $0<\eps<1, \eps_j=\eps,$ раccмотрим конструкцию с параметрами   

$s_j(i)=i$, если   $0< i\leq [(1-\eps)r_j]$,

$s_j(i)=v_j$, если   $[(1-\eps)r_j]< i<r_j$, 

 $s_j(r_j)$ выбираем так, чтобы высоты $h_j$ для $T$ и $S$ были одинаковы.

\vspace{3mm}
Указанная  конструкция $T$ обладает $\eps$-жесткостью: если 
$\|f\|=1, \ \int_Xf d\mu=0, $
 то 
$$(T^{h_j+v_j}f\,,\, f)\to \eps.$$

\vspace{3mm}
 \bf Вынуждение перемешивания. \rm Пусть $\eps_k\to +0$,  фиксированны функции $f_1, f_2, \dots, $.  Мы находим последовательность конструкций $T_k$ и возрастающую последовательность этапов $j_k$    
такую, что при $n>k$  параметры конструкций $T_n$ и $T_k$ одинаковы для всех $j \leq j_k$  при $m_k=h_{j_k}+v_{j_k}$ будет выполнено
$$(T_k^{m_k}f_i,f_i)\ > \ \psi({m_k}), \ i=1,2,\dots,k.\eqno(k)$$

Конструкция    $T$ по определению на  всех этапах $j \leq j_k$ имеет те же параметры, что и конструкция $T_k$, поэтому она будет перемешивающей (напомним, что свойство перемешивания для нее доазывается также как и для обычной лестничной конструкции).   Благодаря произвольному выбору (сколь угодно быстро растущих)  этапов   $j_k$   мы можем обеспечить выполнение неравенства $(k)$ при подстановке $T$ вместо $T_k$ 
(автоморфизмы $T$ и $T_k$ отличаются на множестве такой  малой меры, что  неравенство сохраняется).  Это завершает доказательство леммы 2.

\section{ Эргодические средние и мягкое перемешивание}

Переформулируем слегка теорему 1  работы \cite{P}.

\vspace{2mm}
\bf Теорема 2. \it Пусть для эргодического автоморфизма $T$, некоторой ненулевой  функциии $f\in L_1$ с нулевым средним и  последовательности $k_n\to\infty$  для почти всех $x$ выполнено
$$\frac 1 {k_n} \sum_{i=1}^{k_n}T^if (x)  =o\left(\frac 1 {k_n}\right).$$
Тогда  для почти всех $x$ имеет место сходимость
 $$\frac 1 {N} \sum_{n=1}^{N}T^{k_n}f(x) \to\ f(x). $$ \rm

\vspace{2mm}
В \cite{P} показано, что перемешивающие автоморфизмы не удовлетворяют условиям теоремы 2, и высказана гипотеза о том, что этим свойством обладают только перемешивающие автоморфизмы.  На самом деле этот класс существенно шире: он включает все мягко перемешивающие автоморфизмы,
т.е. автоморфизмы без  жесткого фактора (термин предложен  в работе \cite{FW}).   Отсутствие жесткого фактора  означает отсутствие непостоянной функции  $\chi\in L_1$ такой, что для некоторой $k_n\to\infty$  выполнено $$\|T^{k_n}\chi -\chi\|_1\to 0.$$

Покажем, что при выполнении условий цитируемой теоремы автоморфизм $T$ не  обладает мягким перемешиванием.  Пусть последовательность функций 
$\frac 1 {N} \sum_{n=1}^{N}T^{k_n}f$ сходятся к $f$ п.в., тогда она сходится   по мере. С учетом абсолютной непрерывности интеграла Лебега  получим сходимость по норме:
$$\left\|\frac 1 {N} \sum_{n=1}^{N}T^{k_n}f -f\right\|_1\to 0, \eqno(2)$$
что является несложным упражнением. Будем считать, что $f$  - вещественная функция (общий случай сводится к этому стандартным способом).  Пусть 
$A=\{x: \ a<f(x)<b\}$, причем выполнено $0<\mu(A)<1$. Так как $f\neq const$, такое множество $A$ найдется. Из $(2)$ вытекает, что 
$$\left\|\frac 1 {N} \sum_{n=1}^{N}T^{k_n}\chi_A -\chi_A\right\|_1\to 0. $$
Но это означает, что для фиксированного  $\eps>0$  для большинства индексов 
$n$ выполнено $\|T^{k_n}\chi_A -\chi_A\|_1<\eps$.  Устремляя $\eps$ к
нулю,  выбираем последовательность $n(i)$ такую, что 
$$\|T^{k_{n(i)}}\chi_A -\chi_A\|_1\to 0, \ k_{n(i)}\to\infty. $$
Таким образом,  $T$ не является мягко перемешивающим.
Замечаем также, что ограничение автоморфизма $T$ на минимальную $T$- инвариантную сигма-алгебру, относительно которой измерима функция $f$, является жестким фактором.  При наличии жесткого фактора быструю сходимость средних обеспечивают кограницы \cite{P}. Подведем итог.

\vspace{2mm}
\bf Теорема 3. \it Эргодический  автоморфизм $T$ удовлетворяет условию теоремы 2 только в том случае, когда он обладает нетривиальным жестким фактором. \rm

\vspace{2mm}

Среди многообразия мягко перемешивающих автоморфизмов, не обладающих
свойством перемешивания,  классическим примером является  преобразование  Чакона \cite{J}.  Среди недавних    примеров отметим пуассоновские надстройки над преобразованиями ранга один с полиномиальными слабыми пределами, предложенными в  \cite{19}.
Частично жесткие автоморфизмы из \S 4 также перемешивают  мягко, но  доказательство  этого факта требует значительных усилий.
 
\vspace{2mm}
\bf Благодарности. \rm  Автор признателен рецензенту за замечания.

\newpage 
{V.\,V.~Ryzhikov}

\bf Generic correlations and ergodic averages
for strongly and mildly mixing automorphisms \rm

For a sequence $\psi(n)\to +0$ and a square integrable non-zero function $f$ for the generic  mixing automorphisms $T$  the set  $\{n:\, |(T^nf,f)|>\psi(n)\}$ is infinite. The  mildly mixing automorphisms $T$ do not have а convergence of nonzero averages $\frac 1 {k_n} \sum_{i=1}^{k_n}T^if (x)$ with the rate  $o\left(\frac 1 {k_n}\right)$.

\vspace{2mm}
\it Key words: \rm correlations,  ergodic  averages, generic mixing automorphisms,  partial and  mild mixing, partial rigidity.
\end{document}